\documentclass{article}
\usepackage{graphicx} 
\usepackage{hyperref}
\usepackage{amsmath,amssymb,color,amsthm,enumitem}
\usepackage{comment}
\usepackage{indentfirst}
\usepackage{microtype}
\usepackage{mathtools}
\usepackage[margin=1.2in]{geometry}
\allowdisplaybreaks
\usepackage{url}
\theoremstyle{plain}
\newtheorem{theorem}{Theorem}[section]
\newtheorem{lemma}[theorem]{Lemma}
\newtheorem{corollary}[theorem]{Corollary}
\newtheorem{proposition}[theorem]{Proposition}

\theoremstyle{definition}
\newtheorem{definition}[theorem]{Definition}

\newcommand\ceil[1]{\left\lceil#1\right\rceil}
\newcommand\floor[1]{\left\lfloor#1\right\rfloor}
\renewcommand{\emptyset}{\varnothing}
\newcommand{\p}{\perp}
\newcommand{\m}{\models}

\makeatletter
\newcommand{\labelinthm}[1]{%
    \label{temp#1}
    \protected@write \@auxout {}{\string \newlabel{#1}{{\emph{\ref{temp#1}}}{\thepage}{\emph{\ref{temp#1}}}{temp#1}{}} }%
}
\makeatother
\newcounter{propcounter}

\title{Non-redundant Independence Atoms}
\author{Attila Sali\\HUN-REN Alfr\'ed R\'enyi Institute of Mathematics and \\Department of Computer Science, BUTE \and Jun Yan\thanks{Was supported by
ERC Advanced Grant 883810 during part of this project.}\\ Faculty of Mathematics, University of Waterloo} 

\begin{document}
\maketitle
\begin{abstract}
   Independence atom is a key concept in relational database theory. In the present paper, we prove that the maximum possible size of a non-redundant set of independence atoms on $n$ attributes is asymptotically $(1+o(1))\binom{n}{\floor{n/2}}$. It turns out that this problem is closely related to two extremal set theory problems, the Katona-Tarj\'an question on ``V''-free families of subsets of $[n]$ for the lower bound, and Kleitman's bound on union-free families for the upper bound.
\end{abstract}

\section{Introduction}
Independence and conditional independence are fundamental concepts in areas
as diverse as artificial intelligence, probability theory, social choice theory, and
statistics \cite{Dawid79,Halpern05,Pearl88}. In the present paper we investigate \emph{independence atoms}, which are specializations of the logic concept introduced in \cite{Depandindep} and studied in detail in \cite{KLV}. 

Let $\{A_1,\ldots,A_n\}$ be a set of $n$ \textit{attributes}. Given disjoint $X,Y\subset\{A_1,\ldots,A_n\}$ and an $n$-ary relation $r$ over $\{A_1,\ldots,A_n\}$, we say that $r$ \textit{satisfies} the independence atom $X \p Y$ if for all tuples $t_1 , t_2 \in r$ there is some tuple
$t \in r$ which matches the values of $t_1$ on all coordinates in $X$ and matches the values of $t_2$ on all coordinates in $Y$. More intuitively, if a relation $r$ satisfies $X \p Y$, then
the occurrences of $X$-values are independent of the occurrences of $Y$-values.

As with any other integrity constraint, it is natural to investigate the implication problem for independence atoms. If $\Sigma$ is a collection of independence atoms and $\sigma$ is an independence atom, then $\Sigma$ implies $\sigma$, denoted by $\Sigma\models\sigma$, if any $n$-ary (database) relation $r$ that satisfies every independence atom in $\Sigma$ also satisfies $\sigma$.
The implication problem was axiomatized by Kontinen, Link, and Väänänen \cite{KLV} (see Theorem~\ref{thm:axioms}). 
They also provided a divide-and-conquer algorithm that when given $(\Sigma,\sigma)$, decides in polynomial time whether $\Sigma\models\sigma$. 

Our focus in this paper is the following extremal combinatorial problem raised by Link \cite{LinkPersonal}. What is the maximum possible size of a non-redundant set of independence atoms on $n$ attributes? 
A set $\Sigma$ of independence atoms is \emph{non-redundant} if, for each atom $\sigma \in \Sigma$, $\Sigma\setminus \{\sigma\}\not\models\sigma$. Equivalently, no atom in $\Sigma$ can be derived from the other atoms. Let $f(n)$ be the largest size of a non-redundant set of independence atoms on $n$ attributes. In an undergraduate research paper, Alland, Sali, and Wu~\cite{matcos25} showed that
\[\Theta\left(\frac{2^n}{\sqrt n}\right)=\sum_{i=1}^{n-1}\binom{i}{\ceil{i/2}}\leq f(n)\leq2^{\floor{2n/3}}\binom{n}{\floor{2n/3}}=\Theta\left(\frac{3^n}{\sqrt{n}}\right).\]
Our main result determines $f(n)$ asymptotically.
\begin{theorem}\label{thm:main}
$f(n)=(1+o(1))\binom{n}{\floor{n/2}}$. More precisely, 
\[\left(1+\frac 2n+o\left(\frac 1{n}\right)\right)\binom{n}{\floor{n/2}}\leq f(n)\leq\left(1+\frac{\sqrt{\pi}}{\sqrt{ 2n}}+O\left(\frac 1{n}\right)\right)\binom{n}{\floor{n/2}}.\]
\end{theorem}

We also establish connections between this problem and two other problems in extremal set theory. A family $\mathcal{F}$ of subsets of $[n]$ is called \textit{``V''-free} if there does not exist distinct $A,B,C\in\mathcal{F}$ with $A\subset B$ and $A\subset C$. Let $g(n)$ be the maximum size of a ``V''-free family of subsets of $[n]$, and let $g_2(n)$ be the maximum size of a ``V''-free family of subsets of $[n]$ that contains only sets of size $\floor{n/2}$ and $\floor{n/2}+1$. Katona and Tarj\'an~ \cite{katona2006extremal} proved that $g_2(n)\leq g(n)\leq(1+2/n)\binom{n}{\floor{n/2}}$, and also provided a construction originating from Graham and Sloane's work on constant weight codes \cite{GrahamSloane}, which proved that $g_2(n)\geq(1+1/n+o(1/n))\binom{n}{\floor{n/2}}$. 

A family $\mathcal{F}$ of subsets of $[n]$ is called \textit{union-free} if there does not exist distinct $A,B,C\in\mathcal{F}$ with $A=B\cup C$. Let $h(n)$ be the maximum size of a union-free family of subsets of $[n]$.  Kleitman~\cite{MR409199} proved that $h(n)\leq(1+O(1/\sqrt n))\binom{n}{\floor{n/2}}$. While it is clear that $h(n)\geq g(n)$, as far as we can tell it is not known whether the correct second-order term for $h(n)$ should be $\Theta(1/n)$, $\Theta(1/\sqrt n)$, or something in between. 

We will prove that
\[g_2(n)\leq f(n)\leq h(n)+O\left(\binom{n}{\ceil{n/3}}\right),\]
so lower bounds on $g_2(n)$ and upper bounds on $h(n)$ imply corresponding bounds on $f(n)$ as well. Like for $h(n)$, it remains open what the asymptotically correct second-order term of $f(n)$ should be.

The rest of the paper is organized as follows. In Section~\ref{sec:background}, we provide more background on relational database and independence atoms, and record Kontinen, Link, and Väänänen's \cite{KLV} axiomatization of the implication problem, as well as their algorithm for deciding whether $\Sigma\models\sigma$. In Section~\ref{sec:upper}, we use the algorithmic characterization of implication problem and the LYM inequality to prove the upper bound for $f(n)$. We also establish the connection between $f(n)$ and $h(n)$, which then gives an asymptotically equivalent upper bound using Kleitman's bound on $h(n)$. In Section~\ref{sec:lower}, we provide two lower bounds for $f(n)$. The first one is recursive and improves the lower bound of Alland, Sali, and Wu~\cite{matcos25} by a linear factor in $n$. The second one is better asymptotically, and uses an improved bound on $g_2(n)$ recently obtained by Zhan, Lu, and Zhao~\cite{zhan2026second}.

\section{Database background}\label{sec:background}
In this section, we first summarize basic concepts from the relational model of databases, then introduce the syntax and semantics of independence atoms, and finally summarize the results from \cite{KLV} that are needed in our proofs.

A relation schema is a finite set  $R=\{A_1, A_2,\ldots, A_n\}$ of $n$ different attributes, where for each $A_i$, $dom(A_i)$ is a set of possible values that can be taken for the $i$-th attribute. As an example, the attribute $A_1$ may be eye color, and $dom(A_1)$ is then the set of possible eye colors. A \textit{tuple} over $R$ is an element in $dom(A_1)\times dom(A_2)\times\cdots\times dom(A_n)$, and a \textit{relation} $r$ over $R$ is a finite set of tuples over $R$, or equivalently a subset of $dom(A_1)\times dom(A_2)\times\cdots\times dom(A_n)$. For $X\subset[n]$, let $t(X)\in\prod_{i\in X}dom(A_i)$ denote the \textit{restriction} of the tuple $t$ to $X$, and let $r(X) = \{t(X):t\in r\}\subset\prod_{i\in X}dom(A_i)$ denote the \textit{restriction} of the relation $r$ to $X$.

\begin{definition}
Let $R$ be a relation schema with $n$ attributes. Let $X,Y$ be disjoint subsets of $[n]$. A relation $r$ over $R$ \textit{satisfies} the independence atom $X \p Y$ if for all tuples $t_1, t_2\in r$, there exists a tuple $t \in r$ such that $t(X)=t_1(X)$ and $t(Y)=t_2(Y)$. 
\end{definition}
In other words, if a relation $r$ satisfies $X \p Y$, then the occurrences of $X$-values are independent of the occurrences of $Y$-values, in the sense that $r(X\cup Y ) = r(X) \times r(Y)$.

\begin{definition}
Let $\Sigma$ be a collection of independence atoms on $n$ attributes, and let $\sigma$ be an independence atom on the same set of attributes. We say that $\Sigma$ \textit{implies} $\sigma$, denoted by $\Sigma\models\sigma$, if any relation $r$ over this set of attributes that satisfies every independence atom in $\Sigma$ also satisfies $\sigma$. 
\end{definition}

The implication problem was axiomatized by Kontinen, Link, and Väänänen \cite{KLV} as follows. 

\begin{theorem}[{\cite[Theorem 1]{KLV}}]\label{thm:axioms}
The following four axioms form a sound and complete axiom system for the implication of independence atoms. That is, $\Sigma\models\sigma$ if and only if $\sigma$ can be derived from the atoms in $\Sigma$ using the following four axioms. 
\stepcounter{propcounter}
    \begin{enumerate}[label = \emph{{\textbf{\Alph{propcounter}\arabic{enumi}}}}]
        \item\labelinthm{axiom:1} $\models X \perp \emptyset$ (trivial independence, $\mathcal{T}$). 
        \item\labelinthm{axiom:2} $X \perp Y \models Y \perp X$ (symmetry, $\mathcal{S}$). 
        \item\labelinthm{axiom:3} $X \perp (Y\cup Z) \models X \perp Y$ (decomposition, $\mathcal{D}$). 
        \item\labelinthm{axiom:4} $X \perp Y \wedge (X\cup Y) \perp Z \models X \perp (Y\cup Z)$ (exchange, $\mathcal{E}$). 
    \end{enumerate}
 That is, $\Sigma\models\sigma$ iff $\sigma$ has a finite derivation from $\Sigma$ using the rules above.
\end{theorem}
We will need the following theorem, which is a straightforward corollary of the algorithmic characterization of implications of independence atoms proven in~{\cite[Theorem 2]{KLV}}. For a set of independence atoms $\Sigma$ and an atom $X\p Y$, let $\Sigma[X\cup Y]$ denote the \textit{restriction} of $\Sigma$ to $X\cup Y$, that is, 
\[
\Sigma[X\cup Y]=\{(U\cap (X\cup Y))\p (V\cap (X\cup Y)):U\p V\in\Sigma\}
\]
\begin{theorem}\label{thm:char}
Let $\Sigma$ be a set of  independence atoms and let $\sigma=X\p Y$ be an atom. Then  $\Sigma\models\sigma$ if and only if one of the following holds.
 \begin{enumerate}
     \item $\sigma$ is trivial, i.e., $X=\emptyset$ or $Y=\emptyset$.
     \item $X\p Y\in\Sigma$ or $Y\p X\in\Sigma$.
     \item There exists a non-trivial $U\p V\in \Sigma[X\cup Y]$ satisfying $U\cup V=X\cup Y$ with the following property. Let $Q=U\cap X$, $R=U\cap Y$, $S=V\cap X$, and $T=V\cap Y$, so that $U = Q\cup R$, $V = S\cup T$, $X = Q\cup S$, and $Y = R\cup T$. Then, both $\Sigma[X\cup Y]\models Q\p R$ and  $\Sigma[X\cup Y]\models S\p T$ hold.
 \end{enumerate}
\end{theorem}

\section{Upper Bound}\label{sec:upper}
In this section, we prove the upper bound part of Theorem~\ref{thm:main}. We begin by defining \textit{minimal} non-redundant sets of independence atoms, and studying their properties. 
\begin{definition}\label{def:minimal}
A non-redundant set $\Sigma$ of independence atoms is called \emph{containment minimal} or simply \emph{minimal}, if for every atom $A\p B\in\Sigma$ and all subsets $A'\subset A$ and $B'\subset B$ satisfying $A'\cup B'\not=A\cup B$, we have $\Sigma\setminus\{A\p B\}\m A'\p B'$.
\end{definition}

\begin{lemma}\label{lemma:minimality}
Suppose $\Sigma$ is a non-redundant set of independence atoms. Let $A\p B\in\Sigma$, let $A'\subset A$, and let $B'\subset B$. If $\Sigma\setminus\{A\p B\}\not\m A'\p B'$, then $\Sigma^*:=(\Sigma\setminus\{A\p B\})\cup\{A'\p B'\}$ is non-redundant. 
\end{lemma}
\begin{proof}
Suppose for a contradiction that $\Sigma^*$ is redundant. Then, we have the following two cases. 

\textbf{Case I.} $\Sigma^*\setminus\{A'\p B'\}\m A'\p B'$. This is impossible as $\Sigma^*\setminus\{A'\p B'\}=\Sigma\setminus\{A\p B\}$, and we have $\Sigma\setminus\{A\p B\}\not\m A'\p B'$ from the assumption.

\textbf{Case II.} There exists $X\p Y\in\Sigma^*$ with $X\p Y\not=A'\p B'$, such that $\Sigma^*\setminus\{X\p Y\}\m X\p Y$. Note that $X\p Y\in\Sigma$ as well, and $A\p B\m A'\p B'$ by~\ref{axiom:2} and~\ref{axiom:3}, so by replacing the occurrences of $A'\p B'$ in the derivation of $\Sigma^*\setminus\{X\p Y\}\m X\p Y$ with $A\p B$, we get a derivation of $\Sigma\setminus\{X\p Y\}\m X\p Y$. This contradicts that $\Sigma$ is non-redundant. 
\end{proof}

Note that by repeated applications of Lemma~\ref{lemma:minimality}, for any non-redundant set $\Sigma$ of independence atoms, there is a minimal non-redundant set $\Sigma'$ of the same size. Therefore, when studying $f(n)$, the maximum size of a non-redundant set, we can often assume it is minimal.  

The following result shows that the underlying sets corresponding to independence atoms in a minimal non-redundant set are pairwise distinct.
\begin{lemma}\label{lemma:inj}
Let $\Sigma$ be a minimal non-redundant set of independence atoms on $[n]$. Then, the map $\phi:\Sigma\to\mathcal{P}([n])$ defined by $\phi(A\p B)=A\cup B$ is injective.
\end{lemma}
\begin{proof}

Suppose for a contradiction that there exist distinct $X\p Y$ and $U\p V$ in $\Sigma$, such that $\phi(X\p Y)=X\cup Y=U\cup V=\phi(U\p V)$. Let $Q=X\cap U$, $R=Y\cap U$, $S=X\cap V$, and $T=Y\cap V$. Let $\Sigma'=\Sigma\setminus\{X\p Y\}$, noting that $\Sigma'\m U\p V=(Q\cup R)\p(S\cup T)$. We have two cases.

\textbf{Case I.} At least one of $\Sigma'\not\m Q\p R$ and $\Sigma'\not\m S\p T$ holds. Without loss of generality, assume $\Sigma'\not\m Q\p R$. This contradicts the minimality of $\Sigma$ as $Q\subset X$ and $R\subset Y$. 

\textbf{Case II.} $\Sigma'\m Q\p R$ and $\Sigma'\m S\p T$. Then, combining $S\p T$ and $U\p V=(Q\cup R)\p(S\cup T)$ gives $S\p(Q\cup R\cup T)$ by~\ref{axiom:2} and~\ref{axiom:4}, then~\ref{axiom:3} gives $S\p(R\cup T)$. Similarly, combining $Q\p R$ and $(Q\cup R)\p(S\cup T)$ gives $Q\p(R\cup S\cup T)$ by~\ref{axiom:4}. Finally, combining $Q\p(R\cup S\cup T)$ and $S\p(R\cup T)$ gives $(Q\cup S)\p(R\cup T)=X\p Y$ by~\ref{axiom:2} and~\ref{axiom:4}. Therefore, $\Sigma'\m X\p Y$, contradicting that $\Sigma$ is non-redundant.
\end{proof}

The next result gives a useful characterization of non-redundancy that will be used in both the upper bound and lower bound proofs of Theorem~\ref{thm:main}. Note that in the reverse direction, $\Sigma$ is not necessarily minimal non-redundant. 
\begin{lemma}\label{lemma:contain}
If $\Sigma$ is a minimal non-redundant set of independence atoms, then the following holds. For all distinct atoms $A\perp B$ and $A'\perp B'$ in $\Sigma$, if $A'\cup B'\subset A\cup B$, then either $A'\cup B'\subset A$ or $A'\cup B'\subset B$. Conversely, if this condition holds and $\Sigma$ contains no trivial atoms, then $\Sigma$ is non-redundant.
\end{lemma}
\begin{proof}
First suppose that $\Sigma$ is a minimal non-redundant set, and both $(A'\cup B')\cap A\ne\emptyset$ and $(A'\cup B')\cap B\ne\emptyset$ hold. Let $Q=A'\cap A$, $R=B'\cap A$, $S=A'\cap B$, and $T=B'\cap B$. Set $U=Q\cup R$ and $V=S\cup T$, so that $A\p B\models U\p V$ by~\ref{axiom:3}. Let $\Sigma'=\Sigma\setminus (A'\p B')$.  Then, $U\p V$ is a non-trivial atom in $\Sigma'[A'\cup B']$. Since $Q,S\subset A'$ and $R,T\subset B'$, the definition of minimality implies that $\Sigma\setminus\{A'\p B'\}\m Q\p R $ and $\Sigma\setminus\{A'\p B'\}\m S\p T $. Thus, Theorem~\ref{thm:char} applies to show that $\Sigma'=\Sigma\setminus\{A'\p B'\}\m A'\p B' $, contradicting the minimality of $\Sigma$. 

On the other hand, assume now that $\Sigma$ contains no trivial atoms, and for all distinct $A\perp B$ and $A'\perp B'$ in $\Sigma$ satisfying $A'\cup B'\subset A\cup B$, we have either $A'\cup B'\subset A$ or $A'\cup B'\subset B$. Suppose for a contradiction that $\Sigma$ is redundant, then there exists a non-trivial $A\p B\in\Sigma$ such that $\Sigma':=\Sigma\setminus\{A\p B\}\m A\p B$. If $B\p A\in\Sigma$, then $B\cup A\subset A\cup B$, but $B\cup A$ is contained in neither $A$ nor $B$, contradicting the assumption. Thus, by Theorem~\ref{thm:char}, we can find a non-trivial $U\p V\in \Sigma'[A\cup B]$ with $U\cup V=A\cup B$, so that with $Q=U\cap A$, $R=U\cap B$, $S=V\cap A$, and $T=V\cap B$, we have both $\Sigma'[A\cup B]\models Q\p R$ and  $\Sigma'[A\cup B]\models S\p T$. However, since $U\p V\in \Sigma'[A\cup B]$, there exists $U'
\p V'\in\Sigma'\subset\Sigma$ such that $U=U'\cap(A\cup B)$ and $V=V'\cap(A\cup B)$. In particular, $U\cup V=A\cup B\subset U'\cup V'$, so from the assumption, either $A\cup B\subset U'$ or $A\cup B\subset V'$. This implies that either $U=\varnothing$ or $V=\varnothing$, contradicting that $U\p V$ is a non-trivial atom in $\Sigma'[A\cup B]$.   
\end{proof}

\begin{theorem}\label{thm:upper}
For all $n\geq6$,
\[f(n)\leq\binom{n}{\floor{n/2}}+\frac1{n+1}\sum_{t=3}^{\ceil{n/2}}\binom{n+1}t=\left(1+\frac{\sqrt{\pi}}{\sqrt{2n}}+O\left(\frac1n\right)\right)\binom{n}{\floor{n/2}}.\]
\end{theorem}
\begin{proof}
Let $\Sigma$ be a maximum non-redundant system on $[n]$. By repeatedly applying Lemma~\ref{lemma:minimality} to $\Sigma$, we may assume that $\Sigma$ is minimal. Then, by Lemma~\ref{lemma:inj}, $\phi:\Sigma\to\mathcal{P}([n])$ given by $\phi(A\p B)=A\cup B$ is injective. In particular, for any $S\in\phi(\Sigma)$, there is a unique atom $A_S\p B_S\in\Sigma$ with $A_S\cup B_S=S$. Let $a_k$ be the number of atoms $A\p B$ in $\Sigma$ with $|A|+|B|=k$.

Fix any $2\leq t\leq n-1$. Let $I_t$ be the set of pairs $(S,T)$ such that $S\in\phi(\Sigma)$, $T\subsetneq S$, $|T|=t$, and $T\cap A_S,T\cap B_S\not=\varnothing$. By Lemma~\ref{lemma:contain}, if $(S,T)\in I_t$, then $T\not\in\phi(\Sigma)$.

On the one hand, for every $S\in\phi(\Sigma)$ with $|S|=s$, if $|A_S|=a$ and $|B_S|=b$, with $a+b=s$, then there are $\sum_{i=1}^{t-1}\binom{a}{i}\binom{b}{t-i}=\binom{s}{t}-\binom{a}{t}-\binom{b}{t}\geq\binom{s}{t}-\binom{s-1}{t}=\binom{s-1}{t-1}$ ways to pick $T$.

On the other hand, for a fixed $T$ with $|T|=t$, if $(S_1,T),(S_2,T)\in I_t$, then we cannot have $S_1\subset S_2$. Indeed, otherwise we have $T\subset S_1\subset S_2$, and so $S_1\cap A_{S_2},S_1\cap B_{S_2}\not=\varnothing$. By Lemma~\ref{lemma:contain}, this implies that $S_1\not\in\phi(\Sigma)$, a contradiction. Therefore, for a fixed $T$, the set of $S$ with $(S,T)\in I_t$ form an antichain in the poset $\{R:T\subset R\}$. Thus, by the LYM inequality, 
\[\sum_{S:(S,T)\in I_t}\frac1{\binom{n-t}{|S|-t}}\leq1.\]

Summing over all $|T|=t$ not in $\phi(\Sigma)$, and combining the two previous paragraphs, we get
\[\sum_{s=t+1}^n\frac{a_s\binom{s-1}{t-1}}{\binom{n-t}{s-t}}=\sum_{s=t+1}^n\sum_{S\in \phi(\Sigma):|S|=s}\frac{\binom{s-1}{t-1}}{\binom{n-t}{s-t}}\leq\sum_{T\not\in\phi(\Sigma):|T|=t}\sum_{S:(S,T)\in I_t}\frac1{\binom{n-t}{|S|-t}}\leq\binom{n}{t}-a_t.\]
After simplification, this gives 
\[b_t:=\sum_{s=t}^n\frac{a_s}{s\binom{n}{s}}\leq\frac1t.\]
Set $b_{n+1}=0$. Note that $a_t=(b_t-b_{t+1})t\binom nt$ for every $2\leq t\leq n$, so
\begin{align*}
|\Sigma|&=\sum_{t=2}^na_t=\sum_{t=2}^n(b_t-b_{t+1})t\binom nt\\
&=2b_2\binom n2+\sum_{t=3}^nb_t\left(t\binom nt-(t-1)\binom{n}{t-1}\right).
\end{align*}
Using $b_t\leq1/t$, and that $t\binom nt-(t-1)\binom{n}{t-1}$ is only positive when $t\leq\ceil{n/2}$, we get
\begin{align*}
|\Sigma|&=2b_2\binom n2+\sum_{t=3}^nb_t\left(t\binom nt-(t-1)\binom{n}{t-1}\right)\\
&\leq \binom n2+\sum_{t=3}^{\ceil{n/2}}\left(\binom nt-\frac{t-1}t\binom{n}{t-1}\right)\\
&=\binom n2+\sum_{t=3}^{\ceil{n/2}}\binom nt-\sum_{t=3}^{\ceil{n/2}}\binom{n}{t-1}+\sum_{t=3}^{\ceil{n/2}}\frac{\binom{n}{t-1}}t\\
&=\binom{n}{\ceil{n/2}}+\sum_{t=3}^{\ceil{n/2}}\frac{\binom{n}{t-1}}t\\
&=\binom{n}{\floor{n/2}}+\frac1{n+1}\sum_{t=3}^{\ceil{n/2}}\binom{n+1}t=\left(1+\frac{\sqrt{\pi}}{\sqrt{2n}}+O\left(\frac1n\right)\right)\binom{n}{\floor{n/2}},
\end{align*}
as required.
\end{proof}

As mentioned in the introduction, Kleitman proved the following upper bound on $h(n)$, the maximum size of a union-free family of subsets of $[n]$. 
\begin{theorem}[\cite{MR409199}]\label{thm:kleitman}
Let $h(n)$ be the size of the largest family $\mathcal{F}\subset[n]$, such that there does not exist distinct $A,B,C\in\mathcal{F}$ satisfying $A=B\cup C$. Then, 
\[h(n)\leq\binom{n}{\floor{n/2}}+\frac{2^n}{n+1}+O\left(\binom{n}{\floor{n/4}}\right).\]
\end{theorem}

The following result shows that $h(n)$ is essentially an upper bound on $f(n)$, and in particular recovers asymptotically the same upper bound as Theorem~\ref{thm:upper}. An improved upper bound on $h(n)$ could also lead to a better upper bound on $f(n)$. 
\begin{proposition}
$f(n)\leq h(n)+O\left(\binom{n}{\ceil{n/3}}\right)$.
\end{proposition}
\begin{proof}
Let $\Sigma$ be a maximum non-redundant set of independence atoms on $[n]$. By repeated applications of Lemma~\ref{lemma:minimality}, we may assume $\Sigma$ is minimal. By Lemma~\ref{lemma:inj}, the map $\phi:\Sigma\to\mathcal{P}([n])$ defined by $\phi(A\p B)=A\cup B$ is injective.

Let $\mathcal{F}=\phi(\Sigma)$, and note that $|\mathcal{F}|=|\Sigma|$ by injectivity. We claim that there are no distinct sets $X,Y,Z$ in $\mathcal{F}$ with $X=Y\cup Z$, and $Y\cap Z\not=\varnothing$. Suppose for a contradiction that such $X,Y,Z$ exist, and suppose $A\p B$ is the atom in $\Sigma$ satisfying $\phi(A\p B)=X$. Then, by Lemma~\ref{lemma:contain}, $Y\subset A$ or $Y\subset B$, and $Z\subset A$ or $Z\subset B$. Without loss of generality, suppose $Y\subset A$. Since $Y\cap Z\not=\varnothing$ and $A\cap B=\varnothing$, we must have $Z\subset A$ as well. However, this implies that $X=Y\cup Z\subset A$, a contradiction. 

Now let $\mathcal{F}'=\mathcal{F}\cap\left(\bigcup_{k=\ceil{n/3}}^{\floor{2n/3}-1}\binom{[n]}{k}\right)$. If there exist distinct $X,Y,Z$ in $\mathcal{F}'$ with $X=Y\cup Z$, then we must have $Y\cap Z=\varnothing$ from above. However, this implies that $2\ceil{n/3}\leq |Y|+|Z|=|X|\leq\floor{2n/3}-1$, a contradiction. Thus, $|\mathcal{F}'|\leq h(n)$, and so
\[f(n)=|\Sigma|=|\mathcal{F}|\leq|\mathcal{F}'|+\sum_{k=0}^{\ceil{n/3}-1}\binom{n}{k}+\sum_{k=\floor{2n/3}}^n\binom{n}{k}\leq h(n)+O\left(\binom{n}{\ceil{n/3}}\right),\]
as required. 
\end{proof}
\section{Lower Bound}\label{sec:lower}
In this section we give two different lower bound constructions. The first one is a recursive partition-type construction as follows.
\begin{theorem}\label{thm:partition}
Let $n=a+b$, with $a\geq b\geq1$. Let $\Sigma_a$ be a non-redundant set on $[a]$, and let $\Sigma_b$ be a non-redundant set on $[a+1,n]$. Let $S=\{A\p B:\emptyset\ne A\subset [a], \emptyset\ne B\subset[a+1,n]\text{ with }|A\cup B|=\ceil{n/2}\}$. Then, the following set is non-redundant. 
\[\Sigma:=\Sigma_a\cup\Sigma_b\cup S.\]
Consequently, 
\begin{equation}\label{eq:partition-bound}
f(n)\geq\max_{\substack{a+b=n\\a\geq b\geq1}}\left\{f(a)+f(b)+\binom{n}{\ceil{n/2}}-\binom{a}{\ceil{n/2}}-\binom{b}{\ceil{n/2}}\right\}.
\end{equation}

\end{theorem}
\begin{proof}
For an atom $X\p Y$ in $\Sigma_a$, the projection of every atom in $\Sigma\setminus\Sigma_a$ onto $X\cup Y$ is trivial, so $X\p Y$ cannot be derived from $\Sigma\setminus\{X\p Y\}$ by Theorem~\ref{thm:char}. Similarly, no atom in $\Sigma_b$ can be derived from other atoms in $\Sigma$.

Now let $A\p B\in S$. For any atom in $\Sigma_a\cup\Sigma_b$, its projection onto $A\cup B$ is trivial. For every atom $A'\p B'\in S\setminus\{A\p B\}$, its projection onto $A\cup B$ is $(A'\cap A)\p(B'\cap B)$. However, since $|A\cup B|=|A'\cup B'|$, $(A'\cap A)\cup(B'\cap B)\not=A\cup B$, so by Theorem~\ref{thm:char}, $A\p B$ cannot be derived from $\Sigma\setminus\{A\p B\}$. Thus, $\Sigma$ is non-redundant.
\end{proof}

We remark that the lower bound in~\cite{matcos25} is the special case when $b=1$ is always chosen in (\ref{eq:partition-bound}), which leads to their lower bound $f(n)\geq\sum_{i=1}^{n-1}\binom{i}{\ceil{i/2}}$. They also conjectured that $f(5)=12$, however this turns out to be incorrect, as Theorem~\ref{thm:partition} implies the following when we use $a=3$ and $b=2$ in (\ref{eq:partition-bound}). 
\begin{corollary}\label{cor:f5}
$f(5)\geq13$.
\end{corollary}
\begin{proof}
$f(5)\geq f(3)+f(2)+\binom{3}{2}\binom{2}{1}+\binom{3}{1}\binom{2}{2}=3+1+6+3=13$. A non-redundant set of atoms is
\[\{1\p2,2\p 3,1\p 3\}\cup\{4\p 5\}\]\[{}\cup\{12\p 4, 13\p 4, 23\p 4, 12\p 5, 13\p 5, 23\p 5\}\cup\{1\p 45,2\p 45,3\p45\}.\qedhere\]
\end{proof}

Define a function $\widetilde{f}(n)$ by $\widetilde{f}(n)=f(n)$ for $n=1,2,3,4$, and let $\widetilde{f}(n)$ satisfy the recurrence~(\ref{eq:partition-bound}) with equality for $n\geq 5$. Then, Theorem~\ref{thm:partition} implies that $f(n)\geq\widetilde{f}(n)$ for every $n$. We used ChatGPT~\cite{correctedaiproof} to solve this recurrence and obtain the following lower bound on $f(n)$. In particular, it turns out that for $n\geq 8$, the optimal splits are $(a,b)=(n-3,3)$ if $n$ is divisible by 4, $(a,b)=(n-1,1)$ if $n$ is congruent to 1 mod 4, $(a,b)=(n-2,2)$ and $(n-3,3)$ if $n$ is congruent to 2 mod 4, and $(a,b)=(n-1,1)$ and $(n-2,2)$ if $n$ is congruent to 3 mod 4.
\begin{proposition}\label{prop:partition-floors}
     Let $\widetilde{f}(n)$ be defined using the recurrence~(\ref{eq:partition-bound}) with equality. Then, 
\[f(n)\geq\widetilde{f}(n)=
\sum_{j=1}^{n-1}
\binom{j}{\lfloor j/2\rfloor}
+\left\lfloor\frac n2\right\rfloor
+\left\lfloor\frac n4\right\rfloor-2,
\qquad n\geq5.\]
Asymptotically, 
\[f(n)\geq\widetilde{f}(n)=\begin{dcases}
    \left(1+\frac2{3n}+O\left(\frac1{n^2}\right)\right)\binom{n}{\floor{n/2}}, &\text{ if }n\text{ is even},\\
    \left(1+\frac4{3n}+O\left(\frac1{n^2}\right)\right)\binom{n}{\floor{n/2}}, &\text{ if }n\text{ is odd}.
\end{dcases}\]
\end{proposition}

Our next result establishes a connection between $f(n)$ and $g_2(n)$, the size of a maximum ``V''-free family of subsets of $[n]$ that contains only sets of size $\floor{n/2}$ and $\floor{n/2}+1$. 
\begin{proposition}\label{prop:gn}
$f(n)\ge g_2(n)$.
\end{proposition}
\begin{proof}
Let $\mathcal{F}'$ be a maximum ``V''-free family of subsets of $[n]$ that only contains sets of size $\floor{n/2}$ and $\floor{n/2}+1$. Let $\mathcal{F}=\{F^c:F\in\mathcal{F}'\}$, then $|\mathcal{F}|=|\mathcal{F}'|$ and $\mathcal{F}$ contains only sets of size $\ceil{n/2}$ and $\ceil{n/2}-1$. Let $\mathcal{F}_1$ be the sets of size $\ceil{n/2}$ in $\mathcal{F}$ and let $\mathcal{A}=\mathcal{F}\setminus\mathcal{F}_1$. Since $\mathcal{F}'$ is ``V''-free, every $F\in \mathcal{F}_1$ has at most one $A\in\mathcal{A}$ such that $A\subset F$.
Let $\Sigma$ consist of the following atoms.
\begin{enumerate}
    \item For every $F\in\mathcal{F}_1$ that has a unique $A\in\mathcal{A}$ such that $A\subset F$, add $A\p(F\setminus A)$ to $\Sigma$. For every $F\in\mathcal{F}_1$ with no $A\in\mathcal{A}$ satisfying $A\subset F$, pick $f\in F$ arbitrarily and add $\{f\}\p(F\setminus\{f\})$ to $\Sigma$.
    \item For every $A\in\mathcal{A}$, pick $a\in A$ arbitrarily, and add $\{a\}\p(A\setminus\{a\})$ to $\Sigma$.
\end{enumerate}
Then, $\Sigma$ is non-redundant by Lemma~\ref{lemma:contain}, so $f(n)\geq|\mathcal{F}|=|\mathcal{F}'|=g_2(n).$
\end{proof}

Katona and Tarj\'an~{\cite[Theorem 6]{katona2006extremal}} studied the maximum size of ``V''-free families and proved that 
\[\left(1+\frac1n+O\left(\frac1{n^2}\right)\right)\binom{n}{\floor{n/2}}\leq g(n)\leq\left(1+\frac2n\right)\binom{n}{\floor{n/2}}.\]
Their lower bound construction for $g(n)$ uses only sets of size $\floor{n/2}$ or $\floor{n/2}+1$, so the same lower bound holds for $g_2(n)$. 
This implies that the lower bound of Proposition~\ref{prop:gn} is better for odd $n$ than the partition-type one. Let us note, that the family of 13 subsets of $[5]$ without a ``$V$'' constructed in \cite{katona2006extremal} corresponds to the non-redundant family of 13
independence atoms of Corollary~\ref{cor:f5}. 

Recently, Zhan, Lu, and Zhao~\cite{zhan2026second} obtained the correct second-order term in for $g(n)$.
\begin{theorem}[{\cite[Theorem 1.4]{zhan2026second}}]
\[g(n)=\left(1+\frac2n+o\left(\frac1{n}\right)\right)\binom{n}{\floor{n/2}}.\]
\end{theorem}
Again, their new lower bound construction for $g(n)$ uses only sets of size $\floor{n/2}$ or $\floor{n/2}+1$, so the same bound also holds for $g_2(n)$. By Proposition~\ref{prop:gn}, this immediately gives the following lower bound on $f(n)$ as a corollary, which for large $n$ is asymptotically better than the bound from the recursive partition-type construction above.
\begin{corollary}\label{cor:lower}
\[f(n)\geq\left(1+\frac2n+o\left(\frac1{n}\right)\right)\binom{n}{\floor{n/2}}.\]
\end{corollary}
Finally, we note that since Katona and Tarj\'an proved that $g_2(n)\leq g(n)\leq(1+2/n)\binom{n}{\floor{n/2}}$, this is essentially the best lower bound on $f(n)$ that can be obtained using Proposition~\ref{prop:gn}.

\section{AI statements}
Proposition~\ref{prop:partition-floors} was obtained using ChatGPT. Their proofs are properly validated and are included in \cite{correctedaiproof}.

\end{document}